\documentstyle{article}

\newtheorem{theorem}{Theorem}[section]
\newtheorem{lemma}[theorem]{Lemma}
\newtheorem{prop}[theorem]{Proposition}
\newtheorem{cor}[theorem]{Corollary}

\title{Point Configurations and Cayley-Menger Varieties}

\author{Ciprian S. Borcea}

\date{}

\begin{document}

\maketitle

\begin{abstract}

\noindent

Equivalence classes of $n$-point configurations in Euclidean,
Hermitian, and quaternionic spaces are related, respectively,
to classical determinantal varieties of symmetric, general, and
skew-symmetric bilinear forms.

\noindent
Cayley-Menger varieties arise in the Euclidean case, and have
relevance for mechanical linkages, polygon spaces and rigidity
theory.

\noindent
Applications include upper bounds for realizations of planar
Laman graphs with prescribed edge-lengths and examples of
special Lagrangians in Calabi-Yau manifolds.
\end{abstract}

\medskip \noindent
{\bf Introduction.} \
We are concerned, initially, with {\em configurations of
$n$ labeled (or ordered) points in the Euclidean space} $R^d$,
up to equivalence under congruence and similarity (rescaling).
We require at least two points to be distinct and denote by $C_n(R^d)$
the resulting {\em configuration space}, made of such equivalence classes.

\medskip \noindent
Cayley-Menger varieties appear when point configurations are looked upon
as encoded in the information given by the squared distances between
any pair of points (up to proportionality).

\medskip \noindent
Cayley expressed the necessary relations between these squared distances
as the vanishing of certain determinants, and Menger found sufficient
conditions for a set of solutions to actually represent the mutual
squared distances of a point configuration. These conditions amount
to sign requirements on determinants of the same kind  \cite{Ber} \cite{Blu}.

\medskip \noindent
Here, we look at the matter not so much in terms of {\em distance geometry},
as envisaged e.g. in \cite{Blu} \cite{C-H}, but rather in terms of
{\bf algebraic geometry}. In fact, our medium will be mostly that of
{\em complex algebraic-geometry}, since we are also interested in certain
{\em complexifications of configuration spaces} (cf. \cite{Bor}).

\medskip \noindent
Thus, we define the Cayley-Menger variety $CM^{d,n}(C)$ as the
Zariski-closure of $imC_n(R^d)\subset P_{{{n}\choose {2}}-1}(C)$.

\medskip \noindent
This approach will be useful with respect to the geometry of the
{\em real points} as well, and we wish to emphasize that the
{\em real slices} of our varieties (or rather the ``realistic'' parts
thereof, which correspond to configurations) are objects of
interest in other areas:
{\bf rigidity theory, robot arm motion planning, mechanical linkages,
or molecular conformations}.

\medskip \noindent
We explore first the {\em planar case:} $d=2$, which is privileged because
of the identification $R^2=C$, and relates (via linear sections of codimension
$n-1$) to planar polygon spaces and their Calabi-Yau complexifications, as
studied in \cite{Bor}. Then we generalize to arbitrary $d\leq n-1$.
\small
Note that $n$ points in some affine space may span a subspace of
dimension at most $(n-1)$.
\normalsize

\medskip \noindent
Indeed, for a better perspective on configuration spaces, the associated
Cayley-Menger varieties $CM^{d,n}(C)$ should be considered as a full
family:

$$ CM^{1,n}(C)\subset CM^{2,n}(C)\subset ... \subset CM^{n-1,n}(C)
= P_{{{n}\choose {2}}-1}(C) $$

\noindent
as resulting from natural inclusions: \
$R\subset R^2\subset ...\subset R^{n-1}$.

\medskip \noindent
The {\em key fact} is that the above series of inclusions can be
identified with the natural {\em stratification by rank} of the projective
space $P_{{{n}\choose {2}}-1}(C)=P(Sym^2(C^{n-1}))$ of {\em symmetric
$(n-1)\times (n-1)$ complex matrices}. Under this identification,
the configuration space of n points in $R^d$ corresponds with
{\em real symmetric matrices of rank at most $d$ and with
non-negative eigenvalues} (up to scalars).

\medskip \noindent
In a certain sense, the most important configuration space is that
of $n$ points on a line $(d=1)$, since it leads to the
{\bf quadratic Veronese embedding} of $P_{n-2}$, which rules the
geometry of our varieties for all higher dimensions $d$.

\medskip \noindent
In fact, there's a larger picture which encompasses the one above as
its $Z_2$-invariant part.

\medskip \noindent
If we replace $R^d$ and its inner product with $C^d$ and its
{\em Hermitian} inner product, we may consider configuration spaces
$C_n(C^d)$ of $n$ points in $C^d$ (with at least two points distinct),
up to equivalence under translations, {\em unitary transformations},
and rescaling.

\medskip \noindent
All considerations and constructions related to $C_n(R^d)$ have their
counterpart for $C_n(C^d)$, and one obtains what might be called
(if a name be needed) Hermite-Gram varieties:

$$ HG^{1,n}(C)\subset HG^{2,n}(C)\subset ... \subset HG^{n-1,n}(C)
= P_{(n-1)^2-1}(C) $$

\medskip \noindent
In any event, they are readily identified with the stratification of
the projective space of $(n-1)\times (n-1)$ complex matrices by rank.

\medskip \noindent
The {\em real structure} to be considered on $HG^{d,n}(C)$ is the one
corresponding to the anti-holomorphic involution $M\mapsto M^*$ on
matrices, where $M^*$ is the adjoint (i.e. conjugate transpose) of $M$.

\medskip \noindent
$HG^{1,n}$ is the image of a {\bf Segre embedding}
$(P_{n-2})^2 \rightarrow P_{(n-1)^2-1}$. and it rules the geometry
for higher dimensions $d$, much in the same way the quadratic
Veronese image $CM^{1,n}$ does with respect to $CM^{\ast,n}$.

\medskip \noindent
Clearly, the former scenario is the $Z_2$-invariant part of the
latter one (for the $Z_2$-action given by transposition on matrices).

\medskip \noindent
At this point, it becomes mandatory (cf. \cite{Arn1} \cite{Arn2})
to investigate the quaternionic (or {\em hyper-Hermitian}) case.
It leads to the stratification of skew-symmetric (or {\em Pfaffian})
forms in $2(n-1)$ variables by rank $2d,\ d=1,...,n-1$.

\medskip \noindent
From the point of view of the {\em self-adjoint matrices}
involved, the triadic series corresponding to {\bf R}, {\bf C},
and {\bf H} have an {\em octonionic} analogue for $d=1,2,3$  and
$n\leq 4$. This goes no further because of Desargues' theorem,
which requires associativity from (projective) dimension three on.

\medskip \noindent
For specific applications, we return to the Euclidean planar case,
or, with orientation, the Hermitian $d=1$ case, and interpret the
degree of the corresponding Cayley-Menger variety as an upper bound
for the number of realizations of generically rigid graphs with
$2n-3$ edges of given length.

\medskip \noindent
We also complement the approach of \cite{Bor} on polygon spaces and
special Lagrangians in Calabi-Yau manifolds.

\medskip \noindent
The last section regroups the main aspects in a tetradic summary.

\medskip \noindent
Our  presentation is ordered as follows:

\medskip \noindent
1. \ Planar configurations

\medskip \noindent
2. \ Cayley-Menger varieties $CM^{2,n}=CM_{2n-4}$

\medskip \noindent
3. \ A realization with symmetric forms

\medskip \noindent
4. \ Cayley coordinates and Gram coordinates

\medskip \noindent
5. \ The cone of positive semi-definite symmetric forms

\medskip \noindent
6. \ A Hermitian analogue

\medskip \noindent
7. \ A quaternionic analogue

\medskip \noindent
8. \ An octonionic enclave

\medskip \noindent
9. \ Mechanical linkages and linear sections

\medskip \noindent
10. Polygon spaces and Calabi-Yau manifolds

\medskip \noindent
11. Summary ($T \epsilon \tau \rho \alpha \lambda o \gamma \iota \alpha$)

\section{Planar configurations}

\medskip \noindent
We start with the collection $F_n(R^2)$ of all possible choices of
$n$ {\em distinct} points $(p_i)_{1\leq i \leq n}$ in $R^2$. The space
$F_n(R^2)$ is naturally identified with the
complement of the ``thick diagonal" in $(R^2)^n$:

$$ F_n(R^2) = (R^2)^n - \bigcup_{i\neq j}
                        \{ p=(p_i)_{1\leq i \leq n} \ | \ p_i=p_j \} $$

\noindent
For various purposes, this space can be variously compactified.
Here, because we want to retain the {\em metric aspect}
of a configuration and consider two of them as equivalent if
one can be turned into the other by a Euclidean displacement (i.e. isometry)
and rescaling (i.e. similarity), we will have to consider the
{\em orbit space} of
$F_n(R^2)$ under the (diagonal action) of the group consisting of all these
Euclidean transformations in the plane.

\medskip \noindent
Equivalence under translations can be eliminated by {\em choosing} the
origin to be one particular point (say $p_1=0$). This choice of
representatives makes the natural(permutation) action of the symmetric
group ${\cal S}_n$ on $n$-point configurations
less manifest, but has other advantages, and we shall see this contrast
reflected again between Cayley coordinates and Gram coordinates (cf. section
5). The choice converts $n$ points into $(n-1)$ {\em vectors}:
$p_j-p_1,\ j=2,...,n$.

\medskip \noindent
Equivalence under similarities is now simply expressed by passing to the
{\em projective space} $P((R^2)^{n-1})=P(R^{2n-2})=P_{2n-3}(R)$, and we need
to exclude only the most degenerate case, namely when all points coincide.

\medskip \noindent
Thus, the
configuration (or moduli) space we are going to investigate
is the orbit space:

$$ C_n(R^2) = P((R^2)^{n-1})/ O(2,R)=P_{2n-3}(R)/O(2,R) $$

\noindent
where $O(2,R)$ stands for the group of orthogonal transformations in $R^2$,
and the action is induced from the diagonal action on $(R^2)^{n-1}$.

\medskip \noindent
Obviously, $O(2,R)$ consists of rotations or rotations followed by
reflection in a line through the origin; in other words:
it has two connected components, each, topologically, a circle $S^1$,
and the quotient $C_n(R^2)$ should be a space of dimension $2n-4$.

\medskip \noindent
Using the identification $R^2=C$, we can give a much more explicit
description of this orbit space $C_n(R^2)$.

\medskip \noindent
With the choice $p_1=0$, a configuration is described by
$(p_2, ..., p_n)$, which we write as $(z_2, ..., z_n)\in C^{n-1}$
when we want to emphasize that we think of the points $p_i=z_i$ as complex
numbers.

\medskip \noindent
Note now that multiplication with non-zero complex numbers in $C=R^2$
means precisely a similarity followed by some rotation in $R^2$, while
{\em conjugation} in $C$, amounts to reflecting in the real (first)
axis.

\medskip \noindent
Thus the passage from $C^{n-1}$ to the complex projective space
$P(C^{n-1})=P_{n-2}(C)$ corresponds precisely to eliminating the
configurations with all points coinciding and accounting for rescalings
and rotations. All that remains to ``factor out" is the equivalence under
reflection in a line, that is under conjugation. This establishes :

\begin{prop}
The configuration space $C_n(R^2)$ can be identified with
the quotient of a complex projective space $P_{n-2}(C)$ by conjugation:

$$ C_n(R^2) = P_{n-2}(C)/conj $$

\end{prop}

\medskip \noindent
Note that the subscript in $P_{n-2}(C)$ indicates {\em complex dimension}
$n-2$, that is {\em real dimension} $2(n-2)$. The fixed points of the
conjugation in $P_{n-2}(C)$, make up precisely
the real projective space $P_{n-2}(R)\subset P_{n-2}(C)$, and clearly
they represent
the equivalence classes of configurations with all points collinear.
Thus, the configuration space $C_n(R^2)$ is a {\em smooth manifold}
of real dimension $2(n-2)=2n-4$ {\em away from collinear configurations},
which make up a ``bad" locus parametrized by $P_{n-2}(R)$.

\medskip \noindent
How ``bad'' this locus is depends on $n$, since {\em locally} along it
the topology is the product of (germs at zero):
$R^{n-2}\times (R^{n-2}/Z_2)$, with $R^{n-2}/Z_2$ denoting the quotient
of $R^{n-2}$ under $x\mapsto -x$.

\medskip \noindent
In particular, for $n=4$, the quotient space remains non-singular, since
$R^2/Z_2=C/Z_2\approx R^2=C$ (by $z\mapsto z^2$).

\medskip \noindent
{\bf Example 1.1:} \ For $n=3$ the space of all triangles in $R^2$, up to
isometry and similarity is a closed hemisphere, the boundary corresponding to
degenerated triangles (with vertices alligned, but not all three confounded).

\medskip \noindent
{\bf Example 1.2:} \ For $n=4$, the configuration space of quadrilaterals
$C_4(R^2)$ will be, as remarked above, a non-singular {\em fourfold}.
Actually,

$$ C_4(R^2)=P_2(C)/conj \approx S^4 $$

\noindent
but the isomorphism with the four dimensional sphere is not that
obvious, although the Betti numbers are clearly those of a sphere
\cite{Kui} \cite{Mass} \cite{Arn2}. We outline an argument in section 5.

\medskip \noindent
A simple way to remain in the smooth realm is to use $P_{n-2}(C)$
for computations and rephrase the proposition as follows:

\medskip \noindent
{\em Identification of conjugate points gives a double covering:}

$$ P_{n-2}(C) \rightarrow C_n(R^2) $$

\noindent
{\em ramified over} $P_{n-2}(R) \subset P_{n-2}(C)$

\section{Cayley-Menger varieties $CM^{2,n}=CM_{2n-4}$}

\medskip \noindent
The fact which leads to what we will call {\em Cayley-Menger
varieties} is that if we retain (up to proportionality)
only the {\em mutual (squared) distances} between the points $(p_i)$
we obtain the same ${{n}\choose{2}}$ projective coordinates for
equivalent configurations. (We use
{\em squared distances} in order to have {\em  (quadratic) polynomial
expressions} and remain in the realm of algebraic geometry.) In
other words, we have a natural map:

$$ C_n(R^2) \rightarrow P_{{{n}\choose{2}} -1} \ \ \ (p_i) \mapsto
(s_{ij}=d_{ij}^2= |p_i-p_j|^2 )_{ij} $$

\noindent
At this point we allow {\em complex coordinates}
in the target space, and envisage $P_{{{n}\choose {2}}-1}(R)$
simply as the real points of the {\em complex} projective space
$P_{{{n}\choose {2}}-1}(C)$. Thus, the image of the configuration space
will be considered as a subset of  $P_{{{n}\choose {2}}-1}(C)$, and we
enter into {\em complex algebraic geometry} with:

\medskip \noindent
{\bf Definition 2.1:} \ The {\em complex Cayley-Menger variety}
$CM^{2,n}(C)=CM_{2n-4}(C)$ is defined as the Zariski-closure of the image
of the configuration space $C_n(R^2)$ in $P_{{{n}\choose {2}}-1}(C)$.

\medskip \noindent
Recall that, for a given subset in a projective space, the Zariski-closure
means the vanishing locus of all homogeneous polynomials which vanish on
the given subset.

\medskip \noindent
{\bf Definition 2.2:} \ The {\em real Cayley-Menger variety}
$CM_{2n-4}(R)$ is the intersection of $CM_{2n-4}(C)$
with the real locus $P_{{{n}\choose {2}}-1}(R)$, in other words: the
fixed locus of {\em conjugation} on $CM_{2n-4}(C)$.

\medskip \noindent
Obviously $im(C_n(R^2))$ is a part of $CM_{2n-4}(R)$,
but {\em only a part}:
there are points in $CM_{2n-4}(R)$ which {\em do not}
correspond with
mutual squared distances of a configuration (see
Corollary 2.6 and section 5 below ).  Menger's {\em inequalities},
expressing sign conditions for various Cayley-Menger determinants,
are one way of distinguishing the {\em ``realistic part''}
$im(C_n(R^2))$ (which is {\em semi-algebraic}) from the full algebraic
real part $CM_{2n-4}(R)$. In section 5 we'll see the general distinction
expressed as that between positive semi-definite symmetric forms
and real symmetric forms of possibly other signatures (and rank
at most $d$)- up to sign.

\medskip \noindent
{\bf Example 2,1:} \ For $n=3$, $CM_2(C)=P_2(C)$, and $CM_2(R)=P_2(R)$
strictly contains the closed disc $im(C_3(R^2))$.

\medskip \noindent
We do now a rank computation
which gives a better idea of the image involved in our definition.
We use the double covering described in section 1:

$$ P_{n-2}(C) \rightarrow C_n(R^2) $$

\noindent
which continues with the map:

$$ C_n(R^2) \rightarrow  P_{{{n}\choose{2}} -1} $$

\medskip \noindent
With complex (homogeneous) coordinates $(z_2:...:z_n)$ for $P_{n-2}(C)$,
the composition:

$$ P_{n-2}(C)  \rightarrow  P_{{{n}\choose{2}} -1} $$

\noindent
reads:

$$ (z_2:...:z_n) \mapsto (z_2 \bar{z}_2 : ...:z_n \bar{z}_n : ... :
(z_i-z_j)(\bar{z}_i - \bar{z}_j) : ...) $$

\begin{prop}
The differential (i.e. tangent map) of
\ $ P_{n-2}(C)  \rightarrow  P_{{{n}\choose{2}} -1} $ has rank
$2(n-2)$ on $P_{n-2}(C) - P_{n-2}(R)$, and rank $(n-2)$ on $P_{n-2}(R)$.
\end{prop}

\medskip \noindent
{\em Proof:} \
We consider the map at the affine level:
$C^{n-1} \rightarrow R^{{n}\choose{2}}$.
Using partial derivatives
$\partial /\partial z_i, \ \partial /\partial \bar{z}_i$
we obtain a matrix of the form:

$$  \left( \begin{array}{llcllcl}
\bar{z}_2 & 0 & ... & 0 & \bar{z}_2-\bar{z}_3 & ... & \bar{z}_2-\bar{z}_n \\
      z_2 & 0 & ... & 0 &       z_2-z_3       & ... &      z_2-z_n \\
       0  & \ & \ast & \ & \               & \ast \ast & \
\end{array} \right) $$

\noindent
where $ \ast $ stands for the corresponding matrix for variables $z_3,...,z_n$.
If we add the first column to each of the last $(n-2)$ colums, and then change
their sign, and move the first in front of them, we obtain:

$$  \left( \begin{array}{lclllcl}
 0 & ... & 0 & \bar{z}_2 & \bar{z}_3 & ... & \bar{z}_n \\
 0 & ... & 0 &        z_2 &       z_3       & ... &   z_n \\
 \ & \ast & \ & \          & \      & \ast \ast \ast & \
\end{array} \right) $$

\noindent
which makes plain that the rank of \ $\ast $ increases by two,
unless all $z_i$'s are collinear, when it increases by one (excepting
the origin). The statement follows by induction. \ \ $\Box$

\medskip \noindent
We'll see now that virtually the same computation allows us to give
a description of $CM_{2n-4}(C)$, even {\em before} we say anything
about defining equations.

\medskip \noindent
Indeed, the above context {\em ``complexifies''} naturally as follows:

\medskip \noindent
$\bullet$ we use coordinates $(u,v)$ on $P_{n-2}(C)\times P_{n-2}(C)$,
with $u=(u_2:...:u_n), v+(v_2:...:v_n)$ labeled in agreement with our
labeling in Proposition 2.1, and consider the map:

$$ P_{n-2}(C)\times P_{n-2}(C) \rightarrow P_{{{n}\choose {2}}-1}(C) $$

$$ (u,v) \mapsto (u_2v_2:...:u_nv_n:...: (u_i-u_j)(v_i-v_j): ...) $$

\noindent
i.e. $u$ takes the role of $z$ and $v$ that of $\bar{z}$, but with
$u$ and $v$ independent complex homogeneous coordinates, and with
a complex target, the map is regular (holomorphic). It is
straightforward to see it's everywhere defined.

\begin{prop}
The (holomorphic) differential of the above morphism of complex
projective manifolds:

$$ P_{n-2}(C)\times P_{n-2}(C) \rightarrow P_{{{n}\choose {2}}-1}(C) $$

\noindent
has rank $2(n-2)$ away from the diagonal in $(P_{n-2}(C))^2$, and rank
$(n-2)$ on this diagonal, i.e. on

$$ im[ \ P_{n-2}(C) \rightarrow (P_{n-2}(C))^2 \ | \ w\mapsto (w,w) \ ] $$

\end{prop}

\medskip \noindent
{\em Proof:} \ As above, with $\partial/\partial u_i, \
\partial/\partial v_i$ \ instead of $\partial/\partial z_i, \
\partial/\partial \bar{z}_i$. \ \ \ $\Box$

\begin{prop}
The complex Cayley-Menger variety $CM_{2n-4}(C)$ is the image of

$$ P_{n-2}(C)\times P_{n-2}(C) \rightarrow P_{{{n}\choose {2}}-1}(C) $$

\noindent
The map:

$$ (P_{n-2}(C))^2 \rightarrow CM_{2n-4}(C)=(P_{n-2})^2/Z_2 $$

\noindent
is a double covering branched along the diagonal $P_{n-2}(C)\subset
(P_{n-2}(C))^2$.

\end{prop}

\medskip \noindent
{\em Proof:} \ $(P_{n-2}(C))^2$ contains the double covering $P_{n-2}(C)$
of the configuration space $C_n(R^2)$ via the embedding:

$$ P_{n-2}(C) \rightarrow (P_{n-2}(C))^2 \ , \ z\mapsto (z,\bar{z}) $$

\noindent
that is, as the fixed ocus of the anti-holomorphic involution of
$(P_{n-2}(C))^2$ \ : \ $(u,v)\mapsto (\bar{v},\bar{u})$.

\medskip \noindent
The map to $P_{{{n}\choose {2}}-1}(C)$ clearly takes $(u,v)$ and $(v,u)$
to the same image, and thus the image of $(P_{n-2}(C))^2$ contains the
image of $C_n(R^2)$ used in the definition of the Cayley-Menger variety.

\medskip \noindent
Referring now to an elementary fact about uniqueness of holomorphic
extensions of real power series, we see that the vanishing of a polynomial
on $imC_n(R^2)$ implies its vanishing on $im(P_{n-2}(C))^2$, because of
vanishing (after composition) on the real points of $(P_{n-2}(C))^2$
under the anti-holomorphic involution $(u,v)\mapsto (\bar{v},\bar{u})$.

\medskip \noindent
The fact that $im(P_{n-2}(C))^2$ is irreducible concludes the
identification:

$$ im(P_{n-2}(C))^2 = CM_{2n-4}(C) $$

\medskip \noindent
The fact that this image is precisely the $Z_2$ quotient of $(P_{n-2}(C))^2$
under the involution $(u,v)\mapsto (v,u)$ follows from the rank
computation which gives a generic immersion, and the $n=3$ case:

$$ (P_1(C))^2\rightarrow (P_1(C))^2/Z_2=P_2(C)=CM_2(C) $$

\noindent
which says that, given:

$$ (s_{12}=u_2v_2:s_{23}=(u_2-u_3)(v_2-v_3): s_{3_1}=u_3v_3)\in P_2(C) $$

\noindent
there are (counting multiplicities) two solutions in $(P_1(C))^2$ (the
intersection of two ${\cal O}(1,1)$ divisors), obviously symmetric
under the $Z_2$ action. \ \ \ $\Box$

\begin{cor}
The complex Cayley-Menger variety
$CM_{2n-4}(C)\subset  P_{{{n}\choose {2}}-1}(C)$ is an irreducible
projective subvariety of complex dmension $2n-4$ and degree

$$ D^{2,n}= deg(CM_{2n-4}(C))=\frac{1}{2} {{2n-4}\choose {n-2}} $$

\noindent
It is swept-out by an $(n-2)$-parameter family of linear subspaces
$P_{n-2}(C)$ corresponding to the two $P_{n-2}(C)$ fibrations of
$(P_{n-2}(C))^2$ (identified under the $Z_2$-action).

\end{cor}

\medskip \noindent
{\em Proof:} \ The degree formula comes from a volume computation
on $(P_{n-2}(C))^2$. Let $h_1$ and $h_2$ denote the hyperplane
classes on the two factors. The hyperplane class on
$P_{{{n}\choose {2}}-1}(C)$ pulls back to $(h_1+h_2)$, and:

$$ D^{2,n}=deg(CM_{2n-4}(C))=\frac{1}{2}\int_{(P_{n-2}(C))^2}
(h_1+h_2)^{2(n-2)}= \frac{1}{2} {{2n-4}\choose {n-2}} $$

\medskip \noindent
A similar computation yields the sectional genus:

\begin{cor}
For $n\geq 4$, a generic linear section of codimension $2n-5$ cuts
$CM_{2n-4}(C)$ along a smooth curve of genus:

$$ g^{2,n}= 1+\frac{n-4}{2}D^{2,n}= 1+\frac{n-4}{4} {{2n-4}\choose {n-2}} $$

\end{cor}

\begin{cor}
The real Cayley-Menger variety
$CM_{2n-4}(R)\subset P_{{{n}\choose {2}}-1}(R)$
is the union of the two $Z_2$-quotients: \
$P_{n-2}(C)/conj=imC_n(R^2)$ and $(P_{n-2}(R))^2/Z_2$,
glued along their common $P_{n-2}(R)$ ramification locus:

$$ CM_{2n-4}(R)=P_{n-2}(C)/Z_2 \bigcup_{P_{n-2}(R)}
(P_{n-2}(R))^2/Z_2 $$

\end{cor}

\medskip \noindent
Thus, the ``threshold'' between the ``realistic'' i.e. configuration part
$P_{n-2}(C)/conj=imC_n(R^2)$, and the ``fake'' remnant
$(P_{n-2}(R))^2/Z_2$, is made of (the image of) collinear configurations.
We'll see in the sequel that these collinear configurations are
geometrically paramount for arbitrary dimension $d\leq n-1$.

\medskip \noindent
{\bf Example 2.1. revisited:} \ For $n=3$, we saw that
$P_1(C)/conj=imC_3(R^2)$ is a closed hemisphere or, equivalently,
a closed disc. On the other hand, $(P_2(R))^2/Z_2$ is a M\"{o}bius
band. Gluing the two along their boundary $S^1=P_1(R)$ yields
$CM_2(R)=P_2(R)$.

\medskip \noindent
{\bf Example 1.2. revisited:} \ For $n=4$, we anticipated the result
that $P_2(C)/conj=imC_4(R^2)$ is the 4-sphere $S^4$. It is somewhat
simpler to obtain: $(P_2(R))^2/Z_2\approx P_4(R)$ \  \cite{Mass}.
 Thus:

$$ CM_4(R) = S^4 \bigcup_{P_2(R)} P_4(R) $$

\noindent
is made topologically of the realistic part $S^4$, and the closed ``fake''
part $P_4(R)$, glued along their common $P_2(R)$. See end of section 5 for
details.

\section{A realization with symmetric forms}

\medskip \noindent
If we observe that the map with image $CM^{2,n}(C)=CM_{2n-4}(C)$:

$$ (P_{n-2}(C))^2 \rightarrow P_{{{n}\choose {2}}-1}(C) $$

\noindent
is given by the linear system of all {\em symmetric} divisors in

$$ |{\cal O}(1,1)|=H^0((P_{n-2}(C))^2, {\cal O}(1,1)) $$

\noindent
we are led to the following intrinsic realization.

\medskip \noindent
Let $V$ be a  vector space (of dimension $n-1$), and $P(V)$
the associated projective space (of dimension $n-2$).
For our purposes $V$ should be a {\em complex vector space
defined over the real field} i.e. $V=V_C=V_R\otimes_R C$, with
$V_R$ a real vector space. Essentially: $V_R=R^{n-1}$ and
$V=V_C=C^{n-1}$. Consider:

$$ (P(V))^2 \rightarrow P(V\otimes V) \cdots \rightarrow
P(Sym(V\otimes V)) $$

\noindent
where the first map is the Segre embedding, given by
$(u,v)\mapsto u\otimes v$, and the second {\em rational map}
corresponds to the linear projection on symmetric tensors
along skew-symmetric tensors:

$$ V\otimes V = Sym(V\otimes V) \oplus \bigwedge^2(V) \rightarrow
Symm(V\otimes V)=(V\otimes V)/(V\wedge V) $$

$$ u\otimes v \mapsto \frac{1}{2}(u\otimes v + v\otimes u) $$

\medskip \noindent
Note, in particular, that the restriction of the Segre map
$(P(V))^2 \rightarrow P(V\otimes V)$ to the diagonal
$P(V)\subset (P(V))^2$ gives the quadratic Veronese embedding

$$ \nu_2: P(V) \rightarrow P(Sym(V\otimes V)) $$

$$ w \mapsto \nu_2(w)= w\otimes w $$

\medskip \noindent
As announced earlier, $\nu_2(P(V))\subset P(Sym(V\otimes V))$,
that is: the {\em image of this quadratic Veronese
embedding} of $P(V)$ for $V$ of dimension $n-1$,
will be
the {\em key object} (with its real, respectively complex points)
in our description of configuration spaces and Cayley-Menger varieties:

$$ imC_n(R^d)\subset CM^{d,n}(R)\subset CM^{d,n}(C) $$

\noindent
for arbitrary $d\leq n-1$.

\medskip \noindent
For the moment, we make explicit the situation pertaining to $d=1, 2$.

\medskip \noindent
Considering the identification $V\otimes V=Hom(V^*,V)$, we may speak
of {\em symmetric}: \ $M=^tM$, and {\em skew-symmetric} transformations:
\ $M=-^tM$. With $V=C^{n-1}$ identified with $V^*=C^{n-1}$ via the
standard bilinear form, we have actually {\em symmetric and
skew-symmetric matrices} with complex entries.

\medskip \noindent
Thus, we let $V\otimes V= Hom(C^{n-1},C^{n-1})= {\cal M}(n-1)$, and
$Sym(V\otimes V)=Sym^2(n-1)$ be our abbreviations for the space of
all $(n-1)\times (n-1)$ matrices with complex entries, respectively
all symmetric matrices. With these standard coordinates, the above
maps read:

$$ (P_{n-2}(C))^2 \rightarrow P({\cal M}(n-1))=P_{n(n-2)(C)} \cdots
\rightarrow P(Sym^2(n-1)) $$

\noindent
with the Segre map corresponding to $(u,v)\mapsto u\cdot ^tv$ \ : the
multiplication of the column vector $u$ with the row vector $^tv$, and
the composition of the two maps becoming:

$$ (u,v) \mapsto \frac{1}{2}(u\cdot ^tv + v\cdot ^tu) $$

\medskip \noindent
Clearly, the first map covers the locus of matrices of rank $\leq 1$,
which meets the skew-symmetric locus only in zero, and thus projectively
maps further onto the locus of {\em symmetric matrices of rank} $\leq 2$.
This gives our realization:

\begin{prop}
The complex Cayley-Menger varieties:

$$ CM^{1,n}(C)\subset CM^{2,n}(C)=CM_{2n-4}(C)\subset
P_{{{n}\choose {2}}-1}(C) $$

\noindent
can be identified by a linear change of coordinates in
$P_{{{n}\choose {2}}-1}(C)$ with the loci:

$$ {\cal R}^1(n-1)\subset {\cal R}^2(n-1)\subset P(Sym^2(n-1)) $$

\noindent
defined by symmetric $(n-1)\times (n-1)$ complex matrices of rank one,
respectively at most two.

\medskip \noindent
The change of coordinates is defined over $R$, hence the corresponding
statement for real points holds true.
\end{prop}

\medskip \noindent
It is fairly transparent now,
what the general result $(d\leq n-1)$ should be.
Thus, Cayley-Menger varieties are recognized as familiar
objects of algebraic geometry: cf. \cite{Gia} \cite{Har} \cite{H-T}
\cite{JLP} \cite{Ful} \cite{GKZ}.
In particular, we have:

\begin{cor}
The projective {\em dual} of the Cayley-Menger variety $CM^{1,n}(C)$
can be identified with the determinantal hypersurface given by all
singular symmetric $(n-1)\times (n-1)$ matrices:

$$ CM^{1,n}(C)^*\approx P\{ A=\ ^tA \ |\ detA=0 \} \subset P(Sym^2(n-1))=
P_{{{n}\choose {2}}-1}(C) $$
\end{cor}

\begin{cor}
The Cayley-Menger variety $CM^{2,n}(C)$ is the {\em secant}
(or chordal) {\em variety} of $CM^{1,n}(C)$.
\end{cor}

\medskip \noindent
{\em Proof:} \ This is immediate from the correspondence with
symmetric matrices (Prop. 3.1). A line through two points of
${\cal R}^1(n-1)$ consists of linear combinations of two symmetric
matrices of rank one, which combinations have rank at most two, and
therefore lie in ${\cal R}^2(n-1)$. Thus, our secant variety
$S({\cal R}^1(n-1)$ is contained in ${\cal R}^2(n-1)$.

\medskip \noindent
The fact that we get all of  ${\cal R}^2(n-1)$ is
clear for {\em real} symmetric matrices (by orthogonal diagonalization)
and follows over $C$ by permanence of algebraic (analytic) identities.
\ \ \ $\Box$

\medskip \noindent
As might be expected, this result generalizes to $CM^{d,n}(C)$, which
is the variety of secant $(d-1)$ planes to
$CM^{1,n}(C)$. This brings forth the central role of $CM^{1,n}$ i.e.
of the quadratic Veronese embedding (of $P_{n-2}$).

\medskip \noindent
Our secant variety $S(CM^{1,n}(C))=CM^{2,n}(C)=CM_{2n-4}(C)$ is
clearly {\em defective}, in the sense that it only doubles the
dimension of $CM^{1,n}(C)$ (without the usual $+1$), and we have
therefore:

\begin{cor}
Let $Tan(CM^{1,n}(C))$ denote the variety of projective tangent
spaces to $CM^{1,n}(C)\subset P_{{{n}\choose {2}}-1}(C)$. Then:

$$ CM^{2,n}(C)=S(CM^{1,n}(C))=Tan(CM^{1,n}(C)) \subset
P_{{{n}\choose {2}}-1}(C) $$

\end{cor}

\section{Cayley coordinates and Gram coordinates}

\medskip \noindent
The change of coordinates which gives the realization of the
Cayley-Menger varieties by symmetric matrices (cf. Prop.3.1)
is a simple passage from Cayley coordinates to Gram coordinates.

\medskip \noindent
We adopt here the general setting for $n$-point configurations
$p_1,...,p_n\in R^d$.

\medskip \noindent
The {\em Cayley coordinates} are our familiar homogeneous
coordinates $s_{ij},\ 1\leq i\ <j\leq n$ for $P_{{{n}\choose {2}}-1}(C)$,
and for an $n$-point configuration in $R^d$ we have:

$$ s_{ij}=s_{ij}(p)=|p_i-p_j|^2=<p_i-p_j,p_i-p_j> $$

\noindent
with $<\ ,\ >$ denoting the usual inner product in $R^d$.

\medskip \noindent
In order to relate {\em Gram coordinates} to configurations, we have to
{\bf choose} one of the points as origin, so that
{\bf with the choice} $p_1=0$ for example, the Gram coordinates
$a_{ij}, \ 2\leq i\leq j\leq n$ corresponding to a configuration
would be:

$$ a_{ij}=a_{ij}(p)=<p_i-p_1,p_j-p_1>=<p_i,p_j> $$

\medskip \noindent
It is important to observe that {\em permuting} the $n$ points
of a configuration amounts, in Cayley coordinates, to a corresponding
permutation, but in Gram coordinates, if the permutation affects
the chosen origin, we do not have a permutation of these coordinates
any more.

\medskip \noindent
With this {\em caveat}, and our choice decided for $p_1=0$,
we proceed to relating the two sets of coordinates for configurations
in $C_n(R^d)$. This is simply the {\em cosine theorem}:

$$ a_{ij}=\frac{1}{2}(s_{1i}+s_{1j}-s_{ij})\ ,\ 2\leq i\leq j\leq n
\ \ \ \ \ \ (C-G) $$

\noindent
where $s_{ij}=0 \ \mbox{for} \ i=j$, that is: $a_{ii}=s_{1i}$.

\medskip \noindent
Normally, we look at the Gram coordinates as arranged in a
{\em symmetric $(n-1)\times (n-1)$ matrix} $A$ with entries
$a_{ij}=a_{ji}$, while the Cayley coordonates are arranged in a
(bordered) {\em symmetric $(n+1)\times (n+1)$ matrix} $S$ with
entry indices running from zero to $n$, and:

$$ s_{kk}=0,\  s_{0i}=s_{i0}=1,\ \mbox{and} \  s_{ij}=s_{ji}
\ \mbox{for} \ 1\leq i < j\leq n $$

\begin{lemma}
Let a Cayley matrix $S$ and a Gram matrix $A$ be related by (C-G).
Then:

$$ rk(S)=2+rk(A)\ , \ \mbox{and} \ det(S)=(-1)^n 2^{n-1} det(A) $$

\end{lemma}

\medskip \noindent
{\em Proof:} \ Subtract  column $(\ast,1)$ in $S$ from columns
$(\ast,j), \ j=2,...,n$, then subtract row $(1,\ast)$ from rows
$(i,\ast), \ i=2,...,n$, to obtain $-2A$ in the lower right corner.
The lemma becomes obvious on this form. \ \ \ $\Box$

\begin{lemma}
Let  $p_1=0,p_2,...,p_n\in R^d$ represent a point in $C_n(R^d)$.
Let $P$ denote the $d\times (n-1)$ matrix with {\em columns}
$p_2,...,p_n$.

\medskip \noindent
Then, the Gram matrix $A(p)$ of the configuration
is given by $A(p)= \ ^tP \cdot P$, and $A(p)$ is consequently
{\em positive semi-definite of rank at most} $d$.

\medskip \noindent
Conversely, if $A$ is a (non-zero) real symmetric $(n-1)\times (n-1)$ matrix
which is positive semi-definite of rank at most $d$, there is a
configuration $p_1=0,p_2,...,p_n\in R^d$ (with at least two distinct points),
such that $A=A(p)$.

\end{lemma}

\medskip \noindent
{\em Proof:} \ The first part amounts to observing that:

$$ a_{ij}(p)=<p_i,p_j>= \ ^tp_j\cdot p_i $$

\noindent
Clearly $ \ ^tP \cdot P$ is positive semi-definite of rank
at most $rk(P)\leq d$.

\medskip \noindent
For the converse, we use an orthogonal diagonalization of $G$: \ for
some orthogonal matrix  $T=\ ^tT^{-1}\in O(n-1,R)$ we obtain a diagonal
matrix
$D=TAT^{-1}$ with at most $d$ non-zero eigenvalues which are positive,
and we may suppose the eigenvalues listed in decreasing order along the
diagonal.

\medskip \noindent
$D$ has a ``square root'' $D^{1/2}$, with: $ D=D^{1/2}D^{1/2}$ (and which
commutes with all matrices commuting with $D$).
$D^{1/2}$ is diagonal, with
positive square roots for the corresponding positive eigenvaues in $D$
as the only non-zero eigenvalues.  Then:

$$ A=T^{-1}DT=^tTD^{1/2}\cdot ^tD^{1/2}T $$

\medskip \noindent
We may retain the first $d$ rows in $D^{1/2}T$ (the remaining rows
being obviously zero), and call this $d\times (n-1)$ matrix $P$. Then
$A=\ ^tP\cdot P$ and $p_1=0$, together with the columns of $P$ give
the required configuration. \ \ \ $\Box$

\medskip \noindent
{\bf Remark:} \ This lemma shows how to retrieve a representative for
the equivalence class of a configuration in $C_n(R^d)$ when given the
mutual squared distances between its points i.e. the Cayley matrix $S$:
one produces the Gram matrix $A$ by (C-G) and finds $P$ as above.
\small This clearly applies to molecular conformations, in which context
it appears as: the EMBED algorithm \cite{C-H} (6.3, pg. 303)
\normalsize

\medskip \noindent
We can establish now our anticipated correspondence in full generality.

\medskip \noindent
Recall that the Cayley-Menger variety $CM^{d,n}(C)$ is defined
as the Zariski-closure of the image of the configuration space
$C_n(R^d)$ in the complex projective space $P_{{{n}\choose {2}}-1}(C)$
by the map $p\mapsto S(p)$.
given by all squared distances between the points.

\begin{theorem}
The linear transformation
$A: P_{{{n}\choose {2}}-1}(C)\rightarrow P_{{{n}\choose {2}}-1}(C)$
defined by  passage  from Cayley coordinates $S$ to Gram
coordinates $A(S)$, that is:

$$ S\mapsto A(S), \ \ a_{ij}(S)=\frac{1}{2}(s_{1i}+s_{1j}-s_{ij})\ ,
\ 2\leq i\leq j\leq n \ \ \  \ (C-G) $$

\noindent
identifies the family of Cayley-Menger varieties:

$$ CM^{1,n}(C)\subset CM^{2,n}(C)\subset ... \subset CM^{n-1,n}(C)
= P_{{{n}\choose {2}}-1}(C) $$

\noindent
with the determinantal varieties in
$P_{{{n}\choose {2}}-1}(C)=P(Sym^2(C^{n-1})$
given by symmetric $(n-1)\times (n-1)$ complex matrices of rank at
most $d,\ d=1,...,n-1$, that is:

$$ {\cal R}^1(n-1)\subset {\cal R}^2(n-1)\subset ... \subset
{\cal R}^{n-1}(n-1)=P_{{{n}\choose {2}}-1}(C) $$

\noindent
with $A(CM^{d,n}(C))={\cal R}^d(n-1)$.

\medskip \noindent
Consequently, $CM^{d,n}(C)$ is the {\em variety $S_{d-1}(CM^{1,n}(C))$
of secant $(d-1)$-planes to} $CM^{1,n}(C)$.

\medskip \noindent
The (complex) dimension of $CM^{d,n}(C)$ is: \ $dn-{{d+1}\choose {2}}-1$,\
and its degree (for $d\leq n-2$) is given by the formula:

$$ D^{d,n}= deg(CM^{d,n}(C))=\prod_{k=o}^{n-d-2}
\frac{{{n-1+k}\choose {n-d-1-k}}}{{{2k+1}\choose {k}}} $$

\end{theorem}

\medskip \noindent
{\em Proof:} \ We know that $A(imC_n(R^d))\subset {\cal R}^d(n-1)$
consists precisely of real symmetric $(n-1)\times (n-1)$ matrices of
rank at most $d$ and with non-negative eigenvalues. We have to show that
any homogeneous polynomial vanishing on such matrices necessarily vanishes
on all symmetric matrices of rank at most $d$.

\medskip \noindent
It will be enough (by permanence of algebraic identities) to prove this
for real symmetric matrices. By orthogonal diagonalization, every such
matrix of rank at most $d$ is a linear combination of (diagonal) rank
one positive semi-definite symmetric matrices (i.e. lies in the $(d-1)$-plane
spanned by $d$ points on $imC_n(R)\subset CM^{1,n}(C)$). The {\em convex hull}
of these $d$ matrices clearly determines (in the projective picture)
points in $imC_n(R^d)$, and a polynomial vanishing on the latter must
vanish on the whole linear span. Thus $A(CM^{d,n}(C)={\cal R}^d(n-1)$.

\medskip \noindent
The argument above also proves that $CM^{d,n}(C)=S_{d-1}(CM^{1,n}(C)$,
since clearly $S_{d-1}({\cal R}^1(n-1))\subset {\cal R}^d(n-1)$, and
the image contains the real points, hence the equality.

\medskip \noindent
The dimension formula will be apparent from the next proposition on
resolving the singularities of ${\cal R}^d(n-1)$. It agrees,
of course, with the ``naive'' count of real parameters for
$C_n(R^d)$:

$$ (n-1)d-dim_R O(d,R)-1=dn-{{d+1}\choose {2}}-1 $$

\noindent
where (with $p_1=0$) we need $(n-1)d$ parameters for $p_2,...,p_n$,
and we factor out orthogonal transformations and rescaling.

\medskip \noindent
The general degree formula is more elaborate, and we refer to
\cite{H-T} \cite{JLP} \cite{Ful}. \small One may verify by induction that
for $d=2$ this yields our $D^{2,n}$ in Corollary 2.4. \ \ \ $\Box$
\normalsize

\medskip \noindent
Henceforth, we freely substitute ${\cal R}^d(n-1)$ for $CM^{d,n}(C)$,
or conversely.

\begin{cor}
The Cayley-Menger variety $CM^{d,n}(C)$
can be defined by the family
of homogeneous polynomials of degree $d+1$ expressing the vanishing
of all $(d+1)$-minors in the Gram matrix.
\end{cor}

\medskip \noindent
{\bf Remark:} \ In view of Lemma 4.1., another possibility would
be to choose as defining equations the vanishing of all
$(d+3)$-minors in the Cayley matrix. These would be homogeneous
polynomials of degree $d+1, d+2$, and $d+3$.

\begin{prop}
The singular locus of $CM^{d,n}(C)$ is $CM^{d-1,n}(C)$.

\medskip \noindent
A resolution of singularities for $CM^{d,n}(C)$ can be presented as
a $P_{{{d+1}\choose {2}}-1}$-bundle over the Grassmann manifold
$G(n-d-1,n-1)$ of codimension $d$ subspaces in $C^{n-1}$.

\medskip \noindent
In particular, all Cayley-Menger varieties are rational.
\end{prop}

\medskip \noindent
{\em Proof:} \ First, we remark that the Grassmann manifold
$G=G(n-d-1,n-1)$ has dimension $d(n-d-1)$, and the proposition
yields the dimension formula in the theorem above.

\medskip \noindent
We consider the incidence variety:

$$ {\cal I}^d(n-1)=\{ (\Lambda, A)\in G\times {\cal R}^d(n-1) \ |
 \ \Lambda\subset KerA \} $$

\noindent
which projects onto the Grassmannian with fibers:

$$ {\cal I}^d_{\Lambda}(n-1)=P(Sym(\Lambda^{\perp}\otimes \Lambda^{\perp}))
\approx P_{{{d+1}\choose {2}}-1}(C) $$

\medskip \noindent
The projection on ${\cal R}^d(n-1)$ gives the resolution of the latter,
with the exceptional divisor projecting onto ${\cal R}^{d-1}(n-1)$.
\ \ \ $\Box$

\medskip \noindent
Corollary 3.2. generalizes to:

\begin{prop}
The projective dual
of \ \  $CM^{d,n}(C)\subset  P_{{{n}\choose {2}}-1}(C)$ \ \
can be naturally identified with $CM^{n-d-1,n}(C)$.
\end{prop}

\medskip \noindent
In terms of symmetric matrices, that is: {\em quadrics}, one uses the pairing

$$ (A,B) \mapsto Tr(AB) $$

\noindent
and the description of the projective tangent space at a non-singular
point $A\in {\cal R}^d(n-1)$ as all quadrics vanishing on $KerA$.
\ \ \ $\Box$

\section{The cone of positive semi-definite symmetric forms}

\medskip \noindent \small
A positive semi-definite symmetric form on $R^{n-1}$ is a bilinear symmetric
form ${\cal A} : R^{n-1}\times R^{n-1} \rightarrow R$ with
${\cal A}(x,x)\geq 0$ for any $x\in R^{n-1}$. Symmetric
$(n-1)\times (n-1)$ real matrices with non-negative eigenvalues are
identified with positive semi-definite symmetric forms by \ ${\cal A}(x,y)=
<Ax,y>$, for the usual inner product $<\ ,\ >$.
\normalsize

\medskip \noindent
We have seen above that, in Gram coordinates, or, in other words,
under the isomorphism $CM^{d,n}\approx {\cal R}^d(n-1)$, the image
of the configuration space $C_n(R^d)$ is made precisely of positive
semi-definite symmetric forms of rank at most $d$.

\medskip \noindent
If we let $d$ run from 1 to $n-1$, we have a sequence of inclusions:

$$ imC_n(R)\subset imC_n(R^2)\subset ... \subset imC_n(R^{n-1}) $$

\noindent
with the last term identified with the (projective image of)
{\em the convex cone of all (non-zero) positive semi-definite
symmetric forms on} $R^{n-1}$.

\begin{prop}
The extremal rays of this cone
correspond with $imC_n(R)$ i.e. with {\em collinear configurations}.
\end{prop}

\medskip \noindent
{\em Proof:} \ Every non-zero positive semi-definite symmetric form
of rank $d$ is the barycenter of a $(d-1)$-simplex of such forms of
rank 1 (as argued in the proof of Theorem 4.3). \ \ \ $\Box$

\medskip \noindent
It will be convenient here to carry on our considerations in the
{\em vector space} of all real symmetric matrices rather than
the associated (real) projective space.

\medskip \noindent
We consider on $Sym^2(R^{n-1)}$ the Lorentzian form:

$$ L(A,B)= Tr(AB)-Tr(A)Tr(B) $$

\medskip \noindent
Indeed, $L$ is symmetric and has signature $(n-2,1)$, since positive
on traceless forms and negative on multiples of the identity.

\begin{prop}
The cone of positive semi-definite symmetric forms on $R^{n-1}$
lies within the negative cone of the Lorentzian form $L$, except
for its extremal rays (collinear configurations), which lie on
the ``light cone'' $L(A,A)=0$.
\end{prop}

\medskip \noindent
{\em Proof:} \ The form $L$ is clearly invariant under the action of
the orthogonal group $O(n-1,R)$ on symmetric forms: \ $A\mapsto TAT^{-1}$,
and the proposition is obvious for diagonal forms. \ \ \ $\Box$

\begin{cor}
The Lorentzian form $L$ determines a hyperbolic metric on the
projective image of the negative cone, which becomes a model
of a hyperbolic $[(n-1)^2-1]$-dimensional space.
This induces a Riemann metric on all {\em smooth} configuration strata:
$imC_n(R^{d+1})-imC_n(R^d), \ d\geq 1$. \ \ \ $\Box$
\end{cor}

\medskip \noindent
{\bf Remark:} The form $L$ and the induced metrics {\em depend}
on the choice of Gram coordinates and are not invariant under
the full ${\cal S}_n$-action on configurations.

\medskip \noindent
{\bf Example 1.2. once more:}
We can present now an argument (cf. \cite{Kui}; see also \cite{Mass}
\cite{Arn2}) for the topological
result announced in Example 1.2:

$$ C_4(R^2)=P_2(C)/conj\approx S^4 $$

\medskip \noindent
We use an affine chart $Tr(A)=1$, and picture $C_4(R^2)=imC_4(R^2)$
as the boundary of the convex hull of
$imC_4(R)\subset S^4=\{ Tr(A^2)=1 \}$. This boundary is topologically a
4-sphere as well.

\medskip \noindent
The other claim (cf. end of section 2), concerning the closure of the
``fake'' part, namely:

$$ (P_2(R))^2/Z_2\approx P_4(R) $$

\noindent
also becomes transparent at this stage. One can use the fact that $CM_4$
is a cubic hypersurface in $P_5$, and, choosing a point in the
interior of the convex hull of $imC_4(R)$ as ``center'' for the 4-sphere
of realistic points, map a pair of antipodal realistic points to the
``fake'' one  given by the third intersection of the ``diameter'' with
$CM_4(R)$. This gives the (closed) ``fake'' part as: $S^4/Z_2=P_4(R)$.
\ \ \ $\Box$

\section{A Hermitian analogue}

\medskip \noindent
As outlined in the introduction, one may consider configuration
spaces $C_n(C^d)$ for equivalence classes of $n$ points in $C^d$,
where equivalence, this time, is under translations, {\em unitary
transformations} and rescaling. Again, we require that at least two
points be distinct.

\medskip \noindent
In this section we let $<\ ,\ >$ denote the standard {\em Hermitian
inner product} of $C^d$:

$$ <z,w>=\sum_{k=1}^d z_k\bar{w}_k $$

\medskip \noindent
An analogue of Cayley coordinates would still record only the
Euclidean information, but a Hermitian Gram matrix will take
account of the symplectic imaginary part. Thus, upon {\bf choosing}
the origin at the first point $p_1=0$ of a configuration,
we put:

$$ \alpha_{ij}(p)=<p_i,p_j> \ , \ \ 2\leq i,j\leq n $$

\noindent
defining a {\em Hermitian matrix} $A(p)=P^* \cdot P$, where $P$
denotes the $d\times (n-1)$ complex matrix with columns $p_2,...,p_n$.

\medskip \noindent
This gives a map:

$$ A: C_n(C^d) \rightarrow P_{(n-1)^2-1}(C) $$

\noindent
and we define $HG^{d,n}(C)$ to be the Zariski-closure
of $imC_n(C^d)$.

\medskip \noindent
$P_{(n-1)^2-1}(C)$ is to be conceived as the projective space
associated to the vector space of $(n-1)\times (n-1)$ complex
matrices on which we have the {\em anti-holomorphic involution}:
\ $M\mapsto M^*$.
The fixed points of this involution (i.e.
the real points for this real structure) are precisely the
Hermitian matrices $H=H^*$. We'll put $P_{(n-1)^2-1}(\ast R)$
for these real points, when we need to emphasise the distinction
from the ordinary real points $P_{(n-1)^2-1}(R)$.

\medskip \noindent
Replacing orthogonal diagonalization with unitary diagonalization,
all the arguments in Theorem 4.6 carry through and give:

\begin{theorem}
The family of projective varieties:

$$ HG^{1,n}(C)\subset HG^{2,n}(C)\subset ... \subset HG^{n-1,n}(C)
= P_{(n-1)^2-1}(C) $$

\noindent
coincides with the stratification by rank of the projective
space of $(n-1)\times (n-1)$ complex matrices, with $HG^{d,n}(C)$
corresponding with matrices of rank at most $d$.

\medskip \noindent
The image of the configuration space $imC_n(C^d)\subset HG^{d,n}(C)$
consists precisely of Hermitian positive semi-definite matrices
of rank at most $d$.

\medskip \noindent
$HG^{d,n}(C)$ is the variety of $(d-1)$-planes secant to $HG^{1,n}(C)$.
The latter space is the image of the Segre embedding:

$$ (P_{n-2}(C))^2 \rightarrow P_{(n-1)^2-1}(C),\ \
(u,v)\mapsto u\otimes v $$

\medskip \noindent
The (complex) dimension of $HG^{d,n}(C)$ is: \ $2dn-d(d+2)-1$\ , and its
degree (for $d\leq n-2$) is given by the formula:

$$ {\cal D}^{d,n}= deg(HG^{d,n}(C))= \prod_{k=o}^{n-d-2}
\frac{{{n-1+k}\choose {d}}}{{{d+k}\choose {k}}} $$

\medskip \noindent
The Cayley-Menger varieties $CM^{d,n}(C)$ can be identified (via the same
choice for Gram coordinates) with the fixed points of the
corresponding varieties $HG^{d,n}(C)$ under the holomorphic
involution given by transposition on matrices:

$$ CM^{d,n}(C)=HG^{d,n}(C)^{Z_2} $$

\end{theorem}

\medskip \noindent
Again, the dimension formula will be apparent from the statement below
on resolution of singularities, and agrees with the ``naive'' count
of real parameters for $C_n(C^d)$:

$$ (n-1)2d -dim_R U(d)-1 = 2d(n-1)-d^2-1=2dn-d(d+2)-1 $$

\noindent
where $U(d)$ stands for the unitary group in $C^d$.

\medskip \noindent
For the degree formula, we refer to \cite{Ful} 14.4.11.

\medskip \noindent
{\bf Remark:} The real structures induced on $CM^{d,n}(C)$ by
conjugation and adjunction are clearly the same.

\begin{prop}
The singular locus of $HG^{d,n}(C)$ is $HG^{d-1,n}(C)$.

\medskip \noindent
A resolution of singularities for $HG^{d,n}(C)$ can be presented as a
$P_{d(n-1)-1}$-bundle over the Grassmann manifold $G(n-d-1,n-1)$
of codimension $d$ subspaces in $C^{n-1}$.

\medskip \noindent
In particular, all varieties $HG^{d,n}(C)$ are rational.
\end{prop}

\begin{prop}
The symmetric bilinear form: \ $(A,B)\mapsto Tr(A\cdot \ ^tB)$
identifies the projective dual of $HG^{d,n}(C)$ with $HG^{n-d-1,n}(C)$.
\end{prop}

\medskip \noindent
These are classical results: cf. \cite{Har} \cite{GKZ}.

\medskip \noindent
If we let $d$ run from 1 to $n-1$, we have a sequence of inclusions:

$$ imC_n(C)\subset imC_n(C^2)\subset ... \subset imC_n(C^{n-1}) $$

\noindent
with the last term identified with the (projective image of)
{\em the convex cone of all (non-zero) positive semi-definite
Hermitian forms on} $C^{n-1}$.

\begin{prop}
The extremal rays of this cone
correspond with $imC_n(C)$ i.e. with {\em $C$-collinear configurations}.
\end{prop}

\medskip \noindent
We may define a Hermitian form on the space of $(n-1)\times (n-1)$
complex matrices (i.e. linear operators on $C^(n-1)$) by:

$$ {\cal L}(A,B) = Tr(AB^*)-Tr(A)Tr(B^*) $$

\medskip \noindent
${\cal L}$ restricts to a real Lorentzian form on Hermitian operators
(and restricts further to $L$ on real symmetric forms).

\medskip \noindent
The analogue of Proposition 5.2 now reads:

\begin{prop}
The cone of positive semi-definite Hermitian forms on $C^{n-1}$,
i.e. $imC_n(C^{n-1})$,
lies within the negative cone of the Lorentzian form ${\cal L}$
(on Hermitian operators), except
for its extremal rays ($C$-collinear configurations), which lie on
the ``light cone'' ${\cal L}(A,A)=0$.
\end{prop}

\section{A quaternionic analogue}

\medskip \noindent
For Hamilton's quaternions $H$, we use the standard description:

$$ x=a+bi+cj+dk, \ \ i^2=j^2=k^2=ijk=-1 $$

\noindent
We put \ $x^*=a-bi-cj-dk$ for the {\em conjugate}, and this gives:

$$ Re(x)=Re(x^*)=\frac{1}{2}(x+x^*)=a $$

$$ |x|^2=|x^*|^2=xx^*=x^*x=a^2+b^2+c^2+d^2 $$

\medskip \noindent
$C$ has an $S^2$ family of {\em embeddings} in $H$, since any
quaternion of square $-1$ can be used to represent $i\in C$.
The {\bf choice} which makes $i\in C$ be $i\in H$, and
takes 1 and $j$ as a $C$-basis (for scalar multiplication on the right),
gives an identification:

$$ H=C^2 \ \ \ x\mapsto (u,v) $$

via the matrix description of left multiplication:

$$ x= \left( \begin{array}{lr}
                       u & -\bar{v} \\
                       v & \bar{u}
              \end{array} \right)     $$

$$ u=a+bi, \ v=c-di, \ x=a+bi+cj+dk=u+jv $$

\noindent
and then $x^*$ is genuinely the adjoint of $x$.

\medskip \noindent
This choice is meant to relate to a {\em right} vector space
structure on $H^d$ (i.e. scalar multiplication by $x\in H$ is
to the right), which allows then the familiar description of
$H$-linear maps \ $H^d\rightarrow H^d$ \  as $d\times d$ matrices
with entries in $H$. The expression over $C\subset H$
(with the 1,$j$ basis in each factor) replaces each quaternionic
entry with its corresponding $2\times 2$ block. Adjunction, that
is transposition and conjugation, is then {\em consistent} in the
$H$-version and $C$-version.

\medskip \noindent
Consider now the {\em hyper-Hermitian} inner product on $H^d$:

$$ <x,y>= \sum_{i=1}^d y_i^*x_i  \ \in H $$

\noindent
The $H$-linear transformations preserving this inner product make-up
a {\em compact Lie group} traditionally denoted $Sp(d)$, but perhaps more
suggestively described as the {\em  group of hyper-unitary transformations}.
We have:

$$ Sp(d)=\{ \  T=(t_{ij})_{1\leq i,j\leq d} \ | \ T^*T=I_d \ \} $$

\noindent
and $dim_R Sp(d)=d(2d+1)$, since the Lie algebra is given by:

$$ sp(d)=\{ \ \Theta \ | \ \Theta^*+\Theta=0 \ \} $$

\noindent
and $3d+4{{d}\choose {2}}= 3d+2d(d-1)=d(2d+1)$.

\medskip \noindent
We may consider now the {\em configuration space} $C_n(H^d)$ which
consists of equivalence classes of $n$ ordered points in $H^d$, modulo
translations, hyper-unitary transformations, and rescaling. As before,
the ``naive'' count of parameters proposes the real dimension:

$$ (n-1)4d - dim_R Sp(d) -1= 4d(n-1) -d(2d+1) -1 $$

\medskip \noindent
Again, by translation, we make the first point in a configuration
become the origin: $p_1=0$, and arrange the column vectors $p_2,...,p_n$
in a $d\times (n-1)$ matrix $P$ with quaternionic entries. We then have
an associated
Gram matrix: $A(p)=P^*\cdot P$, which is self-adjoint,
with columns spanning an $H$-subspace of dimension at most $d$, and
positive semi-definite:

$$ <A(p)x,y>=<x,A(p)y>, \ \ <A(p)x,x>=<Px,Px>\geq 0 $$

\medskip \noindent
The analogy with symmetric (R), and Hermitian (C) matrices
continues in the sense that quaternionic self-adjoint matrices become
diagonal in a suitable hyper-unitary basis, and with this, all previous
considerations have their quaternionic avatar.

\medskip \noindent
In fact, if we look at the {\em quadratic form} $<Ax,x>$
associated to a self-adjoint operator $A$, we see that it takes values
in $R$, and with $H=C^2=R^4$, we have natural inclusions:

$$ hHer(H^{n-1},H^{n-1}) \subset Her(C^{2(n-1)},C{2(n-1)} \subset
Sym(R^{4(n-1)},R^{4(n-1)}) $$

\noindent
where, as emphasised in \cite{Arn2}, the first space corresponds
to real quadratic forms invariant under the action of $S^3=Sp(1)=SU(2)$
\ ($<Axu,xu>=u^*(x^*Ax)u=(x^*Ax)u^*u=<Ax,x>$, for any quaternion $u$ of
norm one), and the second space corresponds to real quadratic forms
invariant under the action of $S^1=U(1)=SO(2,R)$.

\medskip \noindent
Obviously, this gives inclusions:

$$  C_n(H^d) \subset C_{2n-1}(C^{2d}) \subset C_{4n-3}(R^{4d}) \ \ \ \
\ \ \ \ (H-C-R) $$

\noindent
and one can rely on the notion of Cayley-Menger
variety $CM^{4d,4n-3}(C)$, which is the Zariski-closure
of the last term in the projective
space of complex symmetric matrices, in order to obtain corresponding
complexifications and Zariski=closures for the first two terms.

\medskip \noindent
However, we saw in the previous section a more direct way to define
$HG^{d,n}(C)$, and we want to present a similar approach in the
quaternionic case. Since we need some notation for the envisaged
varieties, we propose $PG^{d,n}$, which suggests both Pfaff-Gram
and Pl\"{u}cker-Grassmann. The reason for these asociations will
become apparent presently.

\medskip \noindent
Indeed, one can establish a correspondance between hyper-Hermitian
(i.e. quaternionic self-adjoint) matrices acting on $H^{n-1}$ and
(a certain real slice of) skew-symmetric complex matrices acting on
$C^{2(n-1)}$.

\medskip \noindent
We have describe above the identification $H=C^2$ given by:

$$ x=a+bi+cj+dk=u+jv=(u,v)=(a+bi,c-di) $$

\noindent
and we consider now an {\em order four R-linear transformation}:

$$ \sigma : H=C^2 \rightarrow H=C^2 \ \ x=(u,v)\mapsto
\sigma(x)=(\bar{v},-u) $$

\noindent
When we consider left multiplication by $x$, and look upon $x$ as
a $2\times 2$ matrix, the effect of $\sigma$ can be described as
a rotation by $\pi/2$, followed by change of sign for the first
column.

\begin{lemma}
Let $A=(a_{ij})_{1\leq i,j \leq n-1}$  be a hyper-Hermitian matrix
of rank $d$, that is: \ $a_{ij}^*=a_{ji}$, and the $H$-subspace
generated by the columns of $A$ (with right scalar multiplication)
is $d$-dimensional.

\medskip \noindent
Let $\sigma(A)=(\sigma(a_{ij}))_{ij}$, and consider $\sigma(A)$ as
a complex $2(n-1)\times 2(n-1)$ matrix, corresponding to $H=C^2$,
that is: replace each quaternionic entry by the corresponding $2\times 2$
block.

\medskip \noindent
Then, $\sigma(A)$ is skew-symmetric of rank $2d$.
\end{lemma}

\medskip \noindent
{\em Proof:} \ The fact that $\sigma(A)$ turns out skew-symmetric is
straightforward, and for the rank comparison, we use:

$$ \tilde{A}=\left( \begin{array}{cc}
                                  A & 0 \\
                                  0 & \overline{A}
                    \end{array} \right) \ \ \ \mbox{and} \ \ \
 \tilde{\sigma(A)}=\left( \begin{array}{cc}
                                  \sigma(A) & 0 \\
                                  0 & \overline{\sigma(A)}
                    \end{array} \right) $$

\noindent
where $A$ is also considered as a $2(n-1)\times 2(n-1)$ complex matrix.
Clearly, as complex matrices, the first has rank $4d$, and the second
$2\cdot rk[\sigma(A)]$. Thus, we have to show that the two matrices have
equal rank.

\medskip \noindent
But this is clear when we rotate the second matrix with $\pi/2$
around its center (which is the same as a transposition and a
permutation of rows corresponding to reflecting in the horizontal
mid-axix), and change signs for every other column. This yields
precisely:

$$               \left( \begin{array}{cc}
                                   0 & A \\
                                   \overline{A} & 0
                    \end{array} \right) \ \ \ \ \ \ \ \Box $$

\medskip \noindent
The operation just described shows at the same time how to
introduce an anti-holomorphic involution on the projective space
of complex skew-symmetric matrices,
with fixed point locus exactly the image by $\sigma$ of hyper-Hermitian
matrices. Expressed at the level of the $2\times 2$ blocks indexed
by $ij, \ 1\leq i,j\leq n-1$, this involution amounts to: rotating
each block around its center by $\pi$, changing signs along the second
diagonal in each block, and ending by conjugation.

\medskip \noindent
With this prepared backgroud, we may consider now the configuration
space $C_n(H^d)$
as embedded in the projective space of complex skew-symmetric matrices,
that is:

$$ C_n(H^d) \hookrightarrow P(\wedge^2C^{2(n-1)})=
P_{{{2n-2}\choose {2}} -1}(C), \ \ \ \ p\mapsto \sigma(A(p)) $$

\noindent
where $A(p)$ is the Gram matrix for $p_1=0,p_2,...,p_n\in H^d$.

\medskip \noindent
{\bf Definition 7.1:} \ The variety $PG^{d,n}(C)$ is the
Zariski-closure of

$$ C_n(H^d) \subset P_{{{2n-2}\choose {2}} -1}(C) $$

\noindent
with real points $PG^{d,n}(R)$ defined according to the above
real structure on skew-symmetric matrices.

\medskip \noindent
We have now a string of results, analogous to the orthogonal (R),
and unitary (C) context:

\begin{theorem}
The family of projective varieties:

$$ PG^{1,n}(C)\subset PG^{2,n}(C) \subset ... \subset PG^{n-1,n}(C)=
 P_{{{2n-2}\choose {2}} -1}(C) $$

\noindent
coincides with the stratification by (even) rank of the projective
space of skew-symmetric $2(n-1)\times 2(n-1)$ complex matrices, with
$PG^{d,n}(C)$ corresponding with matrices of rank at most $2d$.

\medskip \noindent
$PG^{d,n}(C)$ is the variety of $(d-1)$-planes secant to $PG^{1,n}(C)$.
The latter space is the image of the Grassmann-Pl\"{u}cker
embedding:

$$ G(2n-4,2n-2) \hookrightarrow  P_{{{2n-2}\choose {2}} -1}(C) $$

\noindent
of complex codimension two subspaces of $C^{2n-2}$ into
$P(\wedge^{2n-4}C^{2n-2})$.

\medskip \noindent
The (complex) dimension of $PG^{d,n}(C)$ is: \ $4d(n-1)-d(2d+1)-1$,
and its degree (for $d\leq n-2$) is given by the formula:

$$ deg(PG^{d,n}(C))=\frac{1}{2^{2n-2d-3}}\prod_{i=0}^{2n-2d-4}
\frac{{{2n-2+i}\choose {2d+2i+1}}}{{{2i+1}\choose {i}}}=
\prod_{1\leq i\leq j\leq 2n-2d-3} \frac{2d+i+j}{i+j} $$

\end{theorem}

\medskip \noindent
The identification $PG^{1,n}(C)=G(2n-4,2n-2)\approx G(2,2n-2)$
can be made more explicit as follows:

\medskip \noindent
$(i)$ \ first, we have an identification of the configuration
space $C_n(H)$ with the {\em quaternonic projective $(n-2)$-space
with respect to left scalar multiplication}, denoted here $HP_{n-2}$;

\medskip \noindent
$(ii)$ \ then, $HP_{n-2}$ can be identified with
the {\em quaternionic Grassmannian $G_H(n-2,n-1)$ of codimension one
$H$-subspaces of $H^{n-1}$ with respect to right scalar multiplication};

\medskip \noindent
$(iii)$ \ then, $G_H(n-2,n-1)\subset G(2n-4,2n-2)$ from our
identification $H=C^2$;

\medskip \noindent
$(iv)$ \ and finally $G(2n-4,2n-2)=PG^{1,n}(C)$ from the correspondence
between skew-symmetric matrices of rank two and their kernels.

\medskip \noindent
Note that the real structure introduced above on skew-symmetric
matrices gives as real points on $G(2n-4,2n-2)$ precisely the
$H$-subspaces, i.e. $G_H(n-2,n-1)$. In this sense,
the {\em complexification} $HP_{n-2}(C)$ of $HP_{n-2}=HP_{n-2}(R)$
is the complex Grassmannian $G(2n-4,2n-2)\approx G(2,2n-2)$.

\medskip \noindent
For the degree formula, see: \cite{G-K} \cite{H-T} \cite{JLP}.

\begin{prop}
The singular locus of $PG^{d,n}(C)$ is $PG^{d-1,n}(C)$.

\medskip \noindent
A resolution of singularities for $PG^{d,n}(C)$ can be presented
as a $P_{{{2d}\choose {2}}-1}$-bundle over the Grassmann
manifold $G(2(n-d-1),2(n-1))$ of codimension $2d$ subspaces in
$C^{2n-2}$.

\medskip \noindent
In particular, all varieties $PG^{d,n}(C)$ are rational.
\end{prop}

\begin{prop}
The symmetric bilinear form:
 \ $(A,B)\mapsto Tr(AB)$ identifies
the projective dual of $PG^{d,n}(C)$ with $PG^{n-d-1,n}(C)$.
\end{prop}

\medskip \noindent
This follows from the description of the projective tangent space
at a non-singular point $S\in PG^{d,n}(C)$ as the projective
subspace of all skew-symmetric matrices $T$ which, as two-forms,
restrict to zero
on $Ker(S)$, that is: \ $\ ^tyTx=0$, for $x,y\in Ker(S)$. \ \ \ \ $\Box$

\medskip \noindent
In relation with the sequence of inclusions of configuration
spaces:

$$ C_n(H)\subset C_n(H^2)\subset ... \subset C_n(H^{n-1}) $$

\noindent
it will be convenient to revert to their description in terms of
Gram matrices i.e. quaternionic self-adjoint operators, and then
the last term is identified with the (projective image of) the
{\em convex cone of all (non-zero) positive semi-definite
hyper-Hermitian forms on} $H^{n-1}$.

\medskip \noindent
Corresponding to propositions 5.1 and 6.4, we have:

\begin{prop}
The extrmal rays of this cone are given by $C_n(H)$, that is:
$H$-collinear configurations.
\end{prop}

\medskip \noindent
Similarly,

$$ {\cal L}(A,B)= Tr(\frac{1}{2}(AB+BA))-Tr(A)Tr(B) $$

\noindent
is well defined and Lorentzian on hyper-hermitian operators, and allows
one to picture non-collinear quaternionic configurations as plunged in
a hyperbolic space of dimension $2n^2-5n+2$.

\section{An octonionic enclave}

\medskip \noindent
In this section we mention the varieties related to octonions, from the
perspective developed for $R, C$, and $H$.
They are, mercifully\footnote{After three long dramas, a short satyr play.}
, restricted to $d=1,2,3$ and $n\leq 4$.

\medskip \noindent
The Graves-Cayley octaves, or octonions {\bf O}, can be described
in terms of quaternions $H^2={\bf O}$ by:

$$ e^2=-1,\ \ \ x=x_1+x_2e, \ \ \ y=y_1+y_2e $$

$$ xy= (x_1y_1-y_2^*x_2) + (x_2y_1^*+y_2x_1)e $$

\noindent
There is a conjugation: \ $\tilde{x}=x_1^*-x_2e$\ , with:

$$ |x|^2=x\tilde{x}=\tilde{X}x=x_1x_1^*+x_2^*x_2=|x_1|^2+|x_2|^2 $$

\noindent
The {\em associator}: \ $[x,y,z]=(xy)z-x(yz)$ \ vanishes whenever two
arguments are equal or conjugate.

\medskip \noindent
We give now a brief description of varieties denoted $OG^{d,n}$ for
$(d,n)=(1,3),(2,3)$ and $(1,4),(2,4),(3,4)$, which represent the
octonionic analogue of previous constructions.

\medskip \noindent
When one considers {\bf O}-Hermitian matrices, the $2\times 2$ case has
a straightforward determinant:

$$ det\left(\begin{array}{cc}
                     \alpha & x \\
                     \tilde{x} & \beta
               \end{array} \right) = \alpha\beta-x\tilde{x} \in R,
\ \ \ \mbox{for} \ \alpha,\beta \in R, \ x\in {\bf O} $$

\noindent
Thus, upon complexification $(\otimes_R C)$, one has:

$$ OG^{1,3}(C)=P\{ \left(\begin{array}{cc}
                                   a & x_c \\
                                 \tilde{x_c} & b
                             \end{array} \right) \ | \ ab=x_c\tilde{x_c}
 \ \} \approx Q_8\subset P_9(C)=OG^{2,3}(C) $$

\medskip \noindent
For $3\times 3$ {\bf O}-Hermitian matrices:

$$ A=\left( \begin{array}{ccc}
                      \alpha & z & y \\
                      \tilde{z} & \beta & x \\
                      \tilde{y} & \tilde{x} & \gamma
                   \end{array} \right) \ \ \ \ \ \alpha,\beta,\gamma
\in R,\ x,y,z \in {\bf O} $$

\noindent one considers the {\em Jordan algebra structure} on the
complexification, with commutative product of matrices:

$$ A_c\cdot B_c=\frac{1}{2}(A_cB_c+B_cA_c) $$

\noindent
and determinanat:

$$ det(A_c)=\frac{1}{3}(TrA_c^3)-\frac{1}{2}(TrA_c)(TrA_c^2)+
\frac{1}{6}(TrA_c)^3 $$

\noindent \small
One recognizes on the right hand side the expression of the product
of three variables (the `eigenvalues') in terms of three basic
symmetric functions: sum, sum of squares, and sum of cubes.
\normalsize

\medskip \noindent
Then:

$$  OG^{2,4}(C)=P\{ A_c\ | \ detA_c=0 \}\subset P_{26}(C)=OG^{3,4}(C) $$

$$ OG^{1,4}(C)=sing(OG^{2,4}(C))\subset OG^{2,4}(C) $$

\noindent
$OG^{1,4}(C)$ has an interpretation as ``rank one'' matrices, and is
better known as the {\em (complexified) octonionic projective plane},
with the {\em rational homogeneous space} description: \ $E_6/P$,
where $E_6$ is the exceptional complex simple Lie group of that type,
and $P$ is the maximal parabolic subgroup corresponding to the first
(or last) root in the $E_6$ graph. \cite{Fre} \cite{Lan}

\medskip \noindent
As expected, $OG^{2,4}(C)$ is, on the one hand, the secant variety
of $OG^{1,4}(C)$, and can be identified, on the other hand, with
the projective dual of $OG^{1,4}(C)$.

\section{Mechanical linkages and linear sections}

\medskip \noindent
It was suggested in the introduction that the formalism of Cayley-Menger
varieties would be useful with respect to {\em mechanical linkages}.
Indeed, mechanical linkages are point configurations with constraints
expressed as prescription of (squared) distances between certain pairs
of points (and visualised as {\em rigid bars} connecting those points).

\medskip \noindent
Cayley coordinates become particularly relevant in this context,
since the {\em configuration space of a linkage with $n$ vertices
and $k$ bars in $R^d$} is simply a
{\em linear section of codimension} $(k-1)$ in $P_{{{n}\choose {2}}-1}(C)$
intersecting $imC_n(R^d)\subset CM^{d,n}(R)\subset CM^{d,n}(C)$
with equations involving {\em only the Cayley coordinates
corresponding to the specified rigid bars}.

\medskip \noindent
We wish to emphasise that, for many purposes, {\em  linkage linear sections}
{\bf may not qualify as generic} among all conceivable linear sections
of a given dimension. In fact, it is rather the peculiar character of such
sections that gives distinctiveness to this study and the related
{\em rigidity theory}.

\medskip \noindent
{\bf Definition 9.1:} \
A {\bf mechanical linkage} $(\Gamma^{d,n}(\sigma), p)$
in a Euclidean space of dimension $d$ is a {\em connected graph}
$\Gamma$ on $n$ vertices (labeled
from 1 to $n$), toghether with an assignment of non-negative numbers
$\{ij\}\mapsto \sigma_{ij}$ for all edges $\{ij\}\in \Gamma$, and
a realization $p$ as a configuration of $n$ labeled
points $p_1,...,p_n\in R^d$ (at least two of them distinct)
with the prescribed squared distances:

$$ |p_i-p_j|^2=<p_i-p_j,p_i-p_j>=\sigma_{ij}\ , \ \mbox{for all edges}\
\{ij\}\in \Gamma $$

\medskip \noindent
We let $|\Gamma|$ stand for the cardinality of $\Gamma$ i.e. the number
of edges in the graph.

\medskip \noindent
Obviously, the set of admissible $\sigma =(\sigma_{ij})_{\{ij\}\in \Gamma}$
may be constrained by the nature of the graph, $d$, and $n$. Note that
we want $p$ to define a point in the configuration space $C_n(R^d)$
where we do not allow all points to be one and the same point. Thus, the
choice $\sigma=0$, will be deemed to have no configuration realization.

\medskip \noindent
{\bf Definiton 9.2:} \
The {\bf configuration space} $C(\Gamma^{d,n}(\sigma))$
of a mechanical linkage $(\Gamma^{d,n}(\sigma), p)$ is the space of all
congruence classes of realizations by $n$ labeled points in $R^d$.

\medskip \noindent
Since one can choose any non-zero bar as unit
of measurement, we see that we have a canonical embedding:

$$ C(\Gamma^{d,n}(\sigma))\subset C_n(R^d) $$

\noindent
and with that, a Zariski-closure in the Cayley-Menger variety
$CM^{d,n}(C)$, which is clearly a linear section:

\begin{prop}
The equations defining the Zariski-closure of the configuration
space $C(\Gamma^{d,n}(\sigma))$
in $CM^{d,n}(C)\subset P_{{{n}\choose {2}}-1}(C)$,
with respect to Cayley coordinates $(s_{ij})$, are simply:

$$  \frac{s_{ij}}{\sigma_{ij}}=\frac{s_{kl}}{\sigma_{kl}}\ , \
\mbox{for all edges}\ \{ij\} , \{kl\} \in \Gamma $$

\end{prop}

\medskip \noindent
Of course, the equations are to be understood as:

$$ \sigma_{kl}s_{ij}=\sigma_{ij}s_{kl}, \ \ \ \ \
 \{ij\} , \{kl\} \in \Gamma $$

\noindent
and the projective subspace $L\Gamma_{\sigma}$ they define in
$P_{{{n}\choose {2}}-1}(C)$
will have codimension $|\Gamma|-1$. \ \ \ $\Box$

\medskip \noindent
{\bf Definition 9.3:} \
The projective variety $L\Gamma_{\sigma}\cap CM^{d,n}(C)$
will be called the {\bf linkage variety} associated to the
mechancal linkage $(\Gamma^{d,n}(\sigma), p)$. It contains, as
its ``realistic'' points, the linkage configuration space
$C(\Gamma^{d,n}(\sigma))$.

\medskip \noindent
It is our contention, to be pursued in \cite{B-S}, that linkage
varieties are instrumental in understanding configuration spaces
of mechanical linkages. Here, we illustrate this point by
obtaining an upper bound for the number of realizations of a
generic planar Laman linkage.

\medskip \noindent
{\bf Definition 9.4:} \
A {\bf planar Laman linkage} is a mechanical linkage
$(\Gamma^{2,n}(\sigma), p)$ in $R^2$, with a (connected) graph
$\Gamma$ on $n$ vertices satisfying:

(i) \ $\Gamma$ has $2n-3$ edges;

(ii) \ for any subset of $k$ vertices, there are at most $2k-3$
edges in $\Gamma$ connecting them.

\medskip \noindent
{\bf Remark:} \ These Laman graphs characterize, in dimension two,
the mechanical linkages which, for generic $\sigma$, are locally
rigid (with a minimum number of bars) \cite{Lam}. Since
$dim_R(CM^{2,n}(R))=2n-4$, it is immediate that one needs at least
$2n-3$ bars in order to obtain isolated points in a generic
linkage configuration space.

\begin{prop}
For generic $\sigma$, the number of
possible realizations (up to congruence) of a planar Laman linkage
$(\Gamma^{2,n}(\sigma),p)$ is bounded by\
$\frac{1}{2}{{2n-4}\choose {n-2}}$.
\end{prop}

\medskip \noindent
{\em Proof:} \ Our bound is the degree of the Cayley-Menger
variety $CM^{2,n}(C)=CM_{2n-4}(C)$ (cf. Corollary 2.4), and the claim
follows form the (refined) B\'{e}zout theorem, considering that, for
generic $\sigma$, all possible realizations are infinitesimally rigid
i.e. isolated not only as points in the configuration space
$C(\Gamma^{2,n}(\sigma))$, but as points in the linkage variety as well.
\ \ \ $\Box$

See also \cite{BS} for more details and the case of arbitrary
dimensions.

\section{Polygon spaces and Calabi-Yau manifolds}

\medskip \noindent
Presently, we are going to use the Cayley-Menger varieties
$CM^{2,n}=CM_{2n-4}$ (and $HG^{1,n}$) in relation to
{\em planar polygonal linkages} and their {\em Calabi-Yau
complexifications}.

\medskip \noindent
The approach here complements the one in \cite{Bor},
which emphasizes toric geometry and expands on matters related
to {\em mirror symmetry}. The latter context motivates the study
of {\em special Lagrangian submanifolds} in Calabi-Yau manifolds,
with particular emphasis on special Lagrangian tori.

\medskip \noindent
For the limited purposes envisaged in this section, we may adopt the
(fairly inclusive) terminology which calls {\bf Calabi-Yau} a
{\bf complex projective manifold with vanishing canonical class
modulo torsion}.

\medskip \noindent
Since all considerations here relate to {\bf real points} of
{\bf Calabi-Yau manifolds defined over R},
we need not dwell on the notion of special Lagrangian submanifold,
normally defined for the finite \'{e}tale coverings where one
actually has
a {\em trivial canonical bundle} i.e. a {\em holomorphic volume form}.
When there's a real structure in this more restricted setting, the real
points do yield a special Lagrangian submanifold - for some adequate
Calabi-Yau metric\footnote{While this seems, in the absence of explicit
Calabi-Yau metrics, an ``easy'' way to get hold of some special
Lagrangians, it faces nevertheless the serious challenge of {\em real
algebraic geometry}: assessing the topology of the real slice.}.

\medskip \noindent
We begin with a simple observation, valid because of the
equality $U(1)=SO(2,R)$ between the unitary group and the
special orthogonal group in the plane $C=R^2$.

\begin{lemma}
There's a natural double covering map:

$$ HG^{1,n}(C) \rightarrow CM^{2,n}(C)=CM_{2n-4}(C) $$

\noindent
ramified over $CM^{1,n}(C)$. and extending the natural double covering:

$$ im C_n(C) \rightarrow imC_n(R^2) $$
\end{lemma}

\medskip \noindent
Of course, this is essentially the double covering described
in Proposition 2.3, considering that $HG^{1,n}(C)$ is the
image of the Segre embedding of $(P_{n-2}(C))^2$, and $CM^{2,n}(C)$
the further projection on symmetric tensors (cf. section 3).
In terms of Gram matrices the map is:

$$ A\mapsto \frac{1}{2}(A+ \ ^tA)\ , \ \ \mbox{for} \ A
\ \mbox{Hermitian of rank one} $$

\medskip \noindent \small
This double covering redresses the loss of orientation
for planar configurations identified by reflection.

\medskip \noindent
Accordingly, planar polygon spaces for unoriented $n$-gons will
be realized in $CM^{2,n}(R)$, while their pull-back to $HG^{1,n}(R)$
by the above double covering will give realizations of the corresponding
spaces for oriented $n$-gons.
\normalsize

\medskip \noindent
Obviously, planar polygon spaces will be
particular cases of linkage configuration spaces, corresponding to
a  {\bf polygonal graph} $\Pi=\Pi^{2,n}$.
We adopt tke standard labeling of edges, namely:

$$ \Pi=\Pi^{2,n}=[\{ 1,2\}, \{ 2,3\},...,\{ i,i+1\},...,
\{ n-1,n\},\{ n,1\}] $$

\medskip \noindent
We have a polygon (configuration) space
$C(\Pi^{2,n}(\sigma))$ once we specify
an admissible {\em edge-length-vector} $q=(q_1,...,q_n)$
and put:

$$ \sigma_{i,i+1}=q_i^2 \ \mbox{for} \ 1\leq i \leq n-1, \ \mbox{and}
\ \sigma_{n,1}=q_n^2 $$

\noindent
for the squared length of the bars. By definition, it consists
(in the unoriented case) of all congruence classes of configurations
$p_1,...,p_n\in R^2$, such that:

$$ |p_j-p_i|^2=\sigma_{ij} \ \ \mbox{for all} \ \{ ij\}\in \Pi $$

\medskip \noindent
Since we may take the perimeter as our scale, we may suppose the
edge-length-vector $q$  standardized by:\ $q_1+...q_n=1$. Then
$q$ is admissible if $0\leq q_i\leq 1/2, \ i=1,...,n$, i.e. no edge
is longer than the sum of the rest.

\medskip \noindent
From the previous section we have:

\begin{prop}
The planar polygon space
$C(\Pi^{2,n}(\sigma))$ can be embedded in the real
Cayley-Menger variety $CM^{2,n}(R)$ as the intersection
of its ``realistic'' part $imC_n(R^2)$ with the codimension $(n-1)$
linear section defined by the equations:

$$ \frac{s_{12}}{\sigma_{12}}=\frac{s_{23}}{\sigma_{23}}=...=
\frac{s_{n-1,n}}{\sigma_{n-1,n}}=\frac{s_{n1}}{\sigma_{n1}}
\ \ \ \ \ \ \ (L\Pi_{\sigma}) $$

\noindent
where $s_{ij}$ are Cayley coordinates in $P_{{{n}\choose {2}}-1}$.
\end{prop}

\medskip \noindent
{\bf Remark:} \ The case $d=2$ has a particularly simple way
of describing the ``realistic'' part $imC_n(R^2)\subset CM^{2,n}(R)$:
it is the part of $CM^{2,n}(R)$ contained in the closure
of the negative cone of the Lorentzian form $L$ (cf. section 5).

\medskip \noindent
As one would expect, the polygon space $C(\Pi^{2,n}(\sigma))$ has
singularities only when the linear section $L\Pi_{\sigma}$ meets
the singular locus $CM^{1,n}(R)\subset CM^{2,n}(R)$, that is, when
the edge-length-vector allows a degeneration of the polygon into a
one-dimensonal configuration. This amounts to a relation of the form:

$$ \sum_{i=1}^n \epsilon_i\cdot q_i=0 \ \ \ \mbox{with}\
\epsilon_i=\pm 1 $$

\medskip \noindent
Thus, when $q$ avoids all ``walls'' of this form,
the polygon configuration space $C(\Pi^{2,n}(\sigma(q)))$
is a smooth $(n-3)$-dimensional manifold, and its topology would
change only when the edge-length-vector parameter $q$
``moves across a wall''.

\medskip \noindent
The fact we want to retain here from \cite{Bor} is that the topology
of $C(\Pi^{2,n}(\sigma))$ can be investigated by separate means:
Morse theory -first and foremost.\footnote{The area function, for
generic $q$, is a Morse function.} This will `pre-empt' the question
about the nature of the real points, when we complexify.

\medskip \noindent
{\bf Example:} \ It is easy to see, intuitively, how to make
$C(\Pi^{2,n}(\sigma))$ into a torus. When one cuts a small
corner of an $(n-1)$-gon, and produces an $n$-gon with the
new edge sufficiently small by comparison with the old edges, the
configuration space will be the product of a circle with the old
configuration space (since the small edge can assume any position
around one end, and the polygon closes-up essentially as for a
null new edge). Thus, one can start with a triangle, perform
a succession of $(n-3)$ such small cuts, and obtain the edges for
a $(S^1)^{n-3}$ polygon configuration space.

\medskip \noindent
The complexification process we are about to consider will offer,
in particular, an illustration for the distinction between linkage
sections (cf. section 9) and more general linear sections.

\medskip \noindent
Indeed, the analysis in \cite{Bor} shows that the polygonal
linkage variety \ $L\Pi_{\sigma}\cap CM^{2,n}(C)$ {\em is
always singular} ( a contraction and $Z_2$ quotient of the
resolved Darboux varieties considered there). However, we
may perturb the linear section (still with real coefficients)
and obtain smooth intersections with $CM^{2,n}(C)$, while
the real locus maintains a connected component diffeomorphic
with the polygon configuration space $C(\Pi^{2,n}(\sigma))$:

\begin{theorem}
Let $q$ be an admissible edge-length-vector away from the walls
described above, and let $\sigma$ stand for the corresponding
squared lengths.

\medskip \noindent
A generic linear section of codimension $(n-1)$, defined over $R$,
and sufficiently close to the linkage section $L\Pi_{\sigma}$,
intersects the Cayley-Menger variety $CM^{2,n}(C)$ along a smooth
$(n-3)$-dimensional Calbi-Yau manifold defined over $R$ whose
real points contain a connected component isomorphic with
the polygon configuration space $C(\Pi^{2,n}(\sigma))$.
\end{theorem}

\medskip \noindent
{\em Proof:} \ Generic linear sections of codimension $(n-1)$
will avoid the $(n-2)$-dimensional singular locus $CM^{1,n}(C)$
and meet $CM^{2,n}(C)$ transversly. This must also be the case
for generic linear sections defined over $R$, since the parameter
locus for singular sections is contained in a proper subvariety
invariant under conjugation.

\medskip \noindent
In order to see that the canonical class (modulo torsion) is
trivial, one looks at the pull-back of the section under the
double covering map (cf. lemma 10.1. above):

$$ HG^{1,n}(C) \rightarrow CM^{2,n}(C) $$

\noindent
But $HG^{1,n}(C)=(P_{n-2}(C))^2$ and the codimension $(n-1)$
section corresponds with a smooth intersection of $(n-1)$
hypersurfaces of bidegree $(1,1)$. \ \ \ $\Box$

\medskip \noindent
{\bf Remark:} \ The case of pentagons $(n=5)$ yields
Enriques surfaces defined over $R$. Renouncing the reality
condition, one obtains the full family of Enriques surfaces
constructed by Reye congruences \cite{Cos}.

\medskip \noindent
We conclude this section with another construction,
in the manner
of \cite{Bry}, of special Lagrangian 3-tori in Calabi-Yau
hypersurfaces of the 4-quadric $G(2,4)$. Again, we look at
real loci in particular cases defined over $R$.

\medskip \noindent
We are going to use {\em not} the usual real structure of
$G(2,4)=Q_4\subset P_5(C)$, but the one explained in Theorem 7.2,
which presents the Grassmannian $G(2,2n-2)=PG^{1,n}(C)$ as the
complexification of the quaternionic projective space
$HP_{n-2}=C_n(H)$.

\medskip \noindent
For $n=3$, we have: \  $HP_1\approx S^4$ \ : the one point
compactification of $H$.

\medskip \noindent
The Calabi-Yau threefolds under consideration are degree four
sections of $G(2,4)=Q_4\subset P_5(C)$. Given that the embedding
of $HP_1$ in the Grassmannian is quadratic, we should look
for 3-tori defined in $H$ by the vanishing of an octic polynomial.

\medskip \noindent
Examples of this kind can be produced as follows: consider the
2-torus $(S^1)^2$ defined in $R^4=H$ by:

$$ x_1^2+x_2^2=1 \ , \ \ x_3^2+x_4^2=1 $$

\noindent
The normal bundle is trivial, hence the circle bundle of
radius $r$ in it defines a 3-torus $(S^1)^3$. For $r < 1$
we have an embedding of this torus in $H$ as:

$$ (a,b,c,d)=(\lambda x_1, \lambda x_2, \mu x_3, \mu x_4) $$

\noindent
with $x_i$ as above and \ $(\lambda -1)^2+(\mu -1)^2 = r^2$.

\medskip \noindent
Elimination yields the octic (non-homogeneous) polynomial equation:

$$ [(a^2+b^2+c^2+d^2)^2 -2r^2(a^2+b^2+c^2+d^2) + (2-r^2)^2]^2 =
8^2(a^2+b^2)(c^2+d^2) $$

\noindent
which vanishes on the intended image. This gives:

\begin{prop}
The family of Calabi-Yau threefolds given by degree four
sections of a smooth quadric $Q_4\subset P_5(C)$ contains
members which allow a real structure with real locus a
3-torus.\ \ \ \ $\Box$
\end{prop}

\section{Summary}

\medskip \noindent
In this summary we uniformize the notation by allowing {\bf K}
to become {\bf R}, {\bf C}, {\bf H}, or {\bf O}, that is: to
designate the real, complex, quaternionic or octonionic numbers.

\medskip \noindent
The common features of these algebraic structures are best
expressed in a theorem of Hurwitz stating that a {\em finite
dimensional real vector space} with:

\noindent
(i) a positive definite inner product,

\noindent
(ii) a (distributive) multiplication with $|xy|=|x|\cdot |y|$, and

\noindent
(iii) a unity

\noindent
must be one of them.

\medskip \noindent
Accordingly, $C_n(K^d)$ will stand for the {\em configuration space}
of $n$ points in $K^d$, and $G_K^{d,n}(C)\subset P_{k(n)}(C)$ for
its {\em Zariski-closure} in a complex projective space
$P_{k(n)}(C)=G_K^{n-1,n}(C)$.
Recall that, for $K=O$, this is only symbolic, and the
pairs $(d,n)$ are restricted to $(1,3),(2,3),(1,4)(2,4),(3,4)$.

\medskip \noindent
Thus, with respect to previous notations:

$$ G_R^{d,n}=CM^{d,n},\ \ G_C^{d,n}=HG^{d,n}, \ \
G_H^{d,n}=PG^{d,n} $$

\noindent
corresponding respectively with {\em determinantal varieties} of
symmetric \ (R), general \ (C), and skew-symmetric (H) forms.

\medskip \noindent
The equivalence relation on configurations is modulo translations,
rescaling, and transformations of $K^d$ given by:

\medskip \noindent
(R) \ the orthogonal group: \ $O(d,R)$, \ with \  $dim_RO(d,R)=
{{d}\choose {2}}$

\medskip \noindent
(C) \ the unitary group: \ $U(d)$, \ with \ $dim_RU(d)=d^2$

\medskip \noindent
(H) \ the hyper-unitary group: \ $Sp(d)$, \  with \ $dim_RSp(d)=
d(2d+1)$.

\medskip \noindent
\begin{tabular}{ccc}   \hline
{\bf K} & $dim_CG_K^{d,n}(C)$ & $dim_CG_K^{n-1,n}(C)=k(n)$
\\ \hline
{\bf R} &  $d(n-1)-{{d}\choose {2}}-1$  & ${{n}\choose {2}}-1$
\\
{\bf C} & $2d(n-1)-d^2-1$ & $(n-1)^2-1=n(n-2)$
\\
{\bf H} & $4d(n-1)-d(2d+1)-1$ & ${{2n-2}\choose {2}}-1=(n-2)(2n-1)$
\\ \hline
\end{tabular}

\medskip \noindent
The projective dual of $G_K^{d,n}(C)$ can be identified with
$G_K^{n-d-1,n}(C)$.

\medskip \noindent
For $d=1$, $C_n(K)=G_K^{1,n}(R)\subset G_K^{1,n}(C)$, and for $n=3,4$
we have:

\medskip \noindent
\begin{tabular}{cllll} \hline
{\bf K} & $C_3(K)$ & $G_K^{1,3}(C)$ & $C_4(K)$ & $G_K^{1,4}(C)$
\\ \hline
{\bf R} & $RP_1=S^1$ & $Q_1\subset P_2(C)$ & $RP_2$ & $P_2(C)\subset P_5(C)$
\\
{\bf C} & $CP_1=S^2$ & $Q_2\subset P_3(C)$ & $CP_2$ &
$(P_2(C))^2\subset P_8(C)$
\\
{\bf H} & $HP_1=S^4$ & $Q_4\subset P_5(C)$ & $HP_2$ &
$G(2,6)\subset P_{14}(C)$
\\
{\bf O} & $OP_1=S^8$ & $Q_8\subset P_9(C)$ & $OP_2$ &
$E_6/P\subset P_{26}(C)$
\\ \hline
\end{tabular}

\medskip \noindent
One recognizes in the last column the {\em four Severi varieties}
\cite{Lan}.

\medskip \noindent
\begin{tabular}{clll} \hline
{\bf K} & $C_n(K)$ & Homogeneous under & $G_K^{1,n}(C)$
\\ \hline
{\bf R} & $RP_{n-2}$ & $O(n-1,R)$ & $P_{n-2}(C)$
\\
{\bf C} & $CP_{n-2}$ & $U(n-1)$ & $P_{n-2}(C)\times P_{n-2}(C)$
\\
{\bf H} & $HP_{n-2}$ & $Sp(n-1)$ & $G(2n-4,2n-2)\approx G(2,2n-2)$
\\ \hline
\end{tabular}

\medskip \noindent
$G_K^{d,n}(C)$ is determined by $G_K^{1,n}(C)$ as its variety of
$(d-1)$ secant planes.

\medskip

\medskip \noindent
{\bf Acknowledgement:} The considerations on planar Laman graphs
are the result of collaboration with Ileana Streinu.

\vspace{0.8in} \noindent
Rider University

\noindent Department of Mathematics

\noindent Lawrenceville, NJ 08648

\noindent U.S.A.

\vspace{0.6in}

\noindent   {\em E-mail address:} \begin{verbatim}
borcea@rider.edu
\end{verbatim}

\end{document}